\documentclass[12pt]{article}

\usepackage{latexsym}
\usepackage{calc}
\usepackage{datetime}
\usepackage{amssymb, amsmath, amsfonts, listings, mathrsfs} 
\usepackage{xcolor} 
\usepackage[normalem]{ulem} 

\usepackage{graphicx}
\usepackage[labelfont=bf]{caption}

\newcounter{hours}\newcounter{minutes}

\DeclareMathOperator*{\argmin}{arg\,min}

\newcommand{\stkout}[1]{\ifmmode\text{\sout{\ensuremath{#1}}}\else\sout{#1}\fi} 

\def\nr{\par \noindent}

\def\Def{\stackrel{\mathrm{def}}{=}}

\def\inter{{\rm int \,}}

\def\dom{{\rm dom \,}}

\def\beq{\begin{equation}}
\def\eeq{\end{equation}}

\def\C{\mathbb{C}}
\def\R{\mathbb{R}}
\def\E{\mathbb{E}}

\def\BI{\begin{itemize}}
\def\EI{\end{itemize}}

\newcommand{\SetEQ}{\setcounter{equation}{0}}
\newcommand{\refLE}[1]{\ensuremath{\stackrel{(\ref{#1})}{\leq}}}
\newcommand{\refEQ}[1]{\ensuremath{\stackrel{(\ref{#1})}{=}}}
\newcommand{\refGE}[1]{\ensuremath{\stackrel{(\ref{#1})}{\geq}}}
\newcommand{\refGT}[1]{\ensuremath{\stackrel{(\ref{#1})}{>}}}

\newtheorem{theorem}{Theorem}
\newtheorem{lemma}{Lemma}
\newtheorem{corollary}{Corollary}

\newtheorem{assumption}{Assumption}
\newtheorem{definition}{Definition}

\newtheorem{example}{Example}
\newtheorem{remark}{Remark}
\newcommand{\proof}{\bf Proof: \rm \nr}
\newcommand{\qed}{\hfill $\Box$ \nr \medskip}
\newcommand{\half}{\mbox{${1 \over 2}$}}

\def\ba{\begin{array}}
\def\ea{\end{array}}
\def\beann{\begin{eqnarray*}}
\def\eeann{\end{eqnarray*}}
\def\bea{\begin{eqnarray}}
\def\eea{\end{eqnarray}}

\def\BT{\begin{theorem}}
\def\ET{\end{theorem}}
\def\BL{\begin{lemma}}
\def\EL{\end{lemma}}
\def\BC{\begin{corollary}}
\def\EC{\end{corollary}}
\def\BE{\begin{example}}
\def\EE{\end{example}}
\def\BD{\begin{definition}}
\def\ED{\end{definition}}
\def\BR{\begin{remark}}
\def\ER{\end{remark}}
\def\BAS{\begin{assumption}}
\def\EAS{\end{assumption}}
\def\BI{\begin{itemize}}
\def\EI{\end{itemize}}

\def\BMP{\begin{minipage}{9.5cm}}
\def\EMP{\end{minipage}}
\def\MPT{\begin{minipage}{11.5cm}}
\def\EPT{\end{minipage}}

\def\la{\langle}
\def\ra{\rangle}

\def\QF{\hspace{5ex} \Box}
\def\QR{\hfill \Box}

\title{
\textbf{Universal Complexity Bounds for Universal Gradient Methods in Nonlinear Optimization} }
\author{
Yurii Nesterov
\thanks{Corvinus Centre for Operations Research at Corvinus Institute for Advanced Studies (Corvinus University of Budapest), and School of Data Sciences (Chinese University of Hong Kong (Shenzhen)).  Professor emeritus at UCLouvain, Belgium. Email: Yurii.Nesterov@uclouvain.be. 
}
}

\date{
September 25, 2025
}

\begin{document}
\maketitle

\abstract{In this paper, we provide the universal first-order methods of Composite Optimization with new complexity analysis. It delivers some universal convergence guarantees, which are not linked directly to any parametric problem class. However, they can be easily transformed into the rates of convergence for the particular problem classes by substituting the corresponding upper estimates for the {\em Global Curvature Bound} of the objective function. We analyze in this way the simple gradient method for nonconvex minimization, gradient methods for convex composite optimization, and their accelerated variant. For them, the only input parameter is the required accuracy of the approximate solution. The accelerated variant of our scheme automatically ensures the best possible rate of convergence simultaneously for all parametric problem classes containing the smooth part of the objective function.}

\vspace{10ex}\noindent

\vspace{3ex}\noindent

\thispagestyle{empty}

\newpage\setcounter{page}{1}

\section{Introduction}
\SetEQ

{\bf Motivation.}
A systematic study of worst-case efficiency bounds for numerical methods of Convex Optimization started from the famous monograph by Nemirovski and Yudin \cite{NY}, where the authors developed a comprehensive Complexity Theory for different parametric classes of optimization problems. Their classification was based on H\"older conditions for the first derivative of the objective function, parameterized by the H\"older degree $\nu \in [0,1]$ and the corresponding H\"older constant:
$$
\ba{rcl}
\| \nabla f(x) - \nabla f(y) \|_* & \leq & L_{\nu} \| x - y \|^{\nu}, \quad x, y \in \R^n
\ea
$$
(notation $f \in \C^{1,\nu}(\R^n)$).

In \cite{NY}, we can find optimal methods for different $\nu \in [0,1]$. However, from the practical point of view, these methods were not perfect since very often they require an auxiliary line-search procedure or the knowledge of some additional information about the objective function.

The first optimal method without any additional requirements and without the line search for the most important practical case $\nu = 1$ was developed in \cite{DAN}. Later, in \cite{NNPaper}, there were presented the optimal methods for $\nu \in (0,1)$ with automatic adjustment to the H\"older degree. However, theses methods need some additional information on the objective function (estimate of the H\"older constant, distance to the minimum). Moreover, for the proper choice of the method's parameters, it is necessary to fix the number of steps of the method in advance.

An importance of automatic adjustment to the best-fitting H\"older power is evident since objective functions can be attributed simultaneously to several H\"older classes. It is almost impossible to decide in advance which of them is more appropriate to our particular objective. 

Hence, several authors proposed different techniques for addressing this issue. In \cite{Lan}, it was shown how this can be done in the framework of Level Method \cite{Level}. However, each iteration of this method is much more difficult than the usual step of gradients schemes.

In \cite{DGN}, it was suggested to treat the H\"older condition in the framework of {\em inexact oracle}. In this approach, it is assumed that the answer of the oracle $(f(x), g(x))$ at point $x$ satisfies the following inequality:
\beq\label{eq-Inexact}
\ba{rcl}
f(y) & \leq & \tilde \delta + f(x) + \la g(x), y - x \ra +  {\tilde L \over 2} \| y - x \|^2, \quad y \in \R^n,
\ea
\eeq
where $\tilde \delta>0$ and $\tilde L>0$ are some inexactness parameters. It was shown that these parameters can be chosen as appropriate functions of $\nu$. Therefore, the function from $\C^{1,\nu}(\R^n)$ can be minimized by gradient methods for $\C^{1,1}(\R^n)$. However, the choice of parameters in these schemes is quite complicated.

Finally, in \cite{UGM}, there were developed "universal" gradient methods for solving the {\em composite miniization problem} \cite{Comp}
$$
\ba{c}
\min\limits_{x \in Q} \; [ \; f(x) + \Psi(x) \;],
\ea
$$
where $\Psi(\cdot)$ is a simple closed convex function. In this paper, the next point $x_+$ has to satisfy the {\em Damped Relaxation Condition}
$$
\ba{rcl}
f(x_+) & \leq & \delta + \min\limits_{y \in Q} \Big[ f(\bar x) + \la g(\bar x), y - \bar x\ra + \half \hat L \| y - \bar x \|^2 \Big],
\ea
$$
where $\delta = {\epsilon \over 2}$ and $\epsilon > 0$ is the required accuracy for the approximate solution of our problem in terms of the function value. In the previous papers \cite{DGN, Lan, NNPaper}, the constant $\delta$ was a complicated function of the class parameters.

After developments in \cite{UGM}, one could have an impression that the problem of constructing the universal first-order methods is solved. Unfortunately, it is not the case. Indeed, all available complexity results are related to the functions from H\"older classes. 
$$
\mbox{\sc
And what about the other problem classes?}
$$
We can imagine that there are many of them. Can we address the performance of our methods on other classes, which are may be even not known today?

On the other hand, there are serious doubts that the whole family of H\"older classes $\C_H(\R^n) \Def \bigcup\limits_{\nu \in [0,1]} \C^{1,\nu}(\R^n)$ represents well the class $\C^1(\R^n)$ of all differentiable functions. Our concerns are supported by the following example.
\BE\label{ex-Cons}
Consider the univariate function $f(x) = f_1(x) + f_2(x)$, where $f_1(x) = \half x^2$ and $f_2(x) = {2 \over 3} |x|^{3/2}$ with $x \in \R$. Note that $f_1 \in \C^{1,1}(\R)$ and $f_2 \in \C^{1/2}(\R)$. On the other hand, for $x = 0$, any $h > 0$, and any $\nu \in [0,1]$, we have
$$
\ba{rcl}
0 \; \leq \; h^{-\nu} [ f'(0+h) - f'(0) ] & = & h^{-\nu} \Big[ h + \sqrt{h}\Big] \; = \; h^{1-\nu} + h^{-\nu + 1/2}.
\ea
$$
Hence, if $\nu \in [0,1)$, then the first term is this sum goes to infinity as $h \to \infty$. If $\nu = 1$, then the second term goes to infinity as $h \to +0$. Therefore, $f \not\in \C_H(\R)$.
\qed
\EE

Thus, we come to the same question again:
\vspace{2ex}

{\sc Can we get complexity bounds for wider classes of objective functions?}

\vspace{2ex}
Our paper gives an affirmative answer to this question. We start from obtaining complexity bounds for some {\em internal characteristic} of our objective function, which has no relations to any standard problem class, and which is defined only by the topology of the function. 
This bound serves as an interface between the method and the particular problem classes. When we decide on the problem class we are interested in, it is enough only to bound the growth of this characteristic using the properties of functions from this class. In this way, we get the complexity estimates for all methods supported by the above analysis.

In our paper, this internal characteristics is the {\em Global Curvature Bound} of the objective function (GCB). We are going to show how all of that work for several {\em universal} optimization schemes, some of which are new.

\vspace{1ex}\noindent
{\bf Contents.} In Section \ref{sc-Mod}, we introduce the notion of the global curvature bound $\hat \mu(\cdot)$ for the general functions. We prove its main properties related to its growth and ability to define a shifted upper quadratic approximation (\ref{eq-Inexact}) for the objective function. Our main statement relates the parameters $\tilde \delta$ and $\tilde L$ in (\ref{eq-Inexact}) with the function $\hat \mu(\cdot)$.

In Section \ref{sc-CompGM}, we explain our approach to the universal analysis of the universal gradient methods. We show how the complexity bounds can be written either in the direct form (using the function $\hat \mu(\cdot)$), or in an inverse form, using the {\em complexity gauge} of the objective function. Both bounds are equivalent. However, their availability depends on the properties of a particular method.

In Section \ref{sc-NonConv}, we present a universal method for non-convex functions aiming to find a point with a small norm of the Gradient Mapping.

In Section \ref{sc-UMSimple}, we present a universal complexity analysis for universal primal and dual gradient methods from \cite{UGM} as applied to convex problem of Composite Minimization.
Their direct complexity guarantees look as follows:
$$
\ba{rcl}
2 \hat \mu \left( 2 \sqrt{D \over k} \right) & \leq & \epsilon,
\ea
$$
where $D$ is the Bregmann distance between the starting point and the solution.
For a particular class $\C^{1,1}(\R^n)$, we have $\hat \mu(t) \leq L_1 t^2$. Substituting this estimate in the above inequality, we get the performance guarantees of our methods for $\C^{1,1}(\R^n)$.

In Section \ref{sc-UFast}, we present a universal complexity analysis for a modified version of Universal Fast Gradient method from \cite{UGM}. Its direct complexity guarantee looks as follows:
$$
\ba{rcl}
(k+1) \hat \mu \left( 2 \left({2 \over k+1} \right)^{3/2} D^{1/2} \right) & \leq & \epsilon.
\ea
$$

Finally, in the last Section \ref{sc-Conc}, we discuss the presented results and the future research directions.

\vspace{1ex}\noindent
{\bf Notation.} In what follows, we work in a finite-dimensional linear space $\E$. Its dual space, the space of linear functions on $\E$, is denoted by $\E^*$. For $x \in \E$ and $g \in \E^*$, we denote by $\la g, x \ra$ the value of linear function $g$ at $x$. For the primal space $\E$, we fix an arbitrary norm $\| \cdot \|$.

\section{Global Curvature Bound}\label{sc-Mod}
\SetEQ

In our approach, we do not fix in advance a particular parametric problem class containing the functional components of our optimization problem. A priory, the only fixed object is the norm $\| \cdot \|$ for measuring the distances in $\E$, and its dual norm
$$
\ba{rcl}
\| g \|_* & = & \max\limits_{x \in \E} \{ \la g, x \ra: \; \| x \| \leq 1 \}, \quad g \in \E^*,
\ea
$$
ensuring Cauchy-Schwartz inequality $|\la g, x \ra| \leq \| x \| \cdot \| g \|_*$.

Let us introduce the notion of {\em Global Curvature Bound} (GCB).
For function $f(\cdot)$ with $\dom f \subseteq \E$, it provides us with a uniform upper bound for the curvature of $f(\cdot)$:
\beq\label{def-USmooth}
\ba{rcl}
\hat \mu_f(t) & = & \sup\limits_{x, y \in \dom f, \atop \alpha \in (0,1)} \Big\{ { | \alpha f(x) + (1-\alpha) f(y) - f(\alpha x + (1-\alpha) y) |  \over \alpha(1-\alpha)} : \; \| x - y \| \leq t \Big\}.
\ea
\eeq
For convex functions, the absolute value in this definition can be dropped.
If function $f(\cdot)$ is linear, then this bound is identically equal to zero. The Global Curvature Bound satisfies the following triangle inequality
$$\ba{rcl}
\hat \mu_{f_1+f_2}(t) & \leq & \hat \mu_{f_1}(t) + \hat \mu_{f_2}(t), \quad t \in \dom \hat \mu_{f_1} \bigcap \dom \hat \mu_{f_2},
\ea
$$ 
which becomes equality when $f_2(\cdot)$ is linear.

If no ambiguity arise, the lower index of function $\hat \mu_f(\cdot)$ is dropped. 
We will prove that GCB of any function $f(\cdot)$
satisfies the following {\em growth condition}:
\beq\label{eq-ModUS}
\ba{rcl}
\hat \mu_f(\beta t) & \geq & \beta^2 \hat \mu_f(t), \quad t \in \dom \hat \mu_f, \; \beta \in [0,1].
\ea
\eeq

Our definition of GCB is similar to the notion of {\em uniform smoothness}, which was introduced in \cite{AP} for convex functions with $\dom f = \E$. However, in \cite{AP}, inequality (\ref{eq-ModUS}) was derived by duality relations with the {\em exact modulus of uniform convexity}
\beq\label{def-UConv}
\ba{rcl}
\mu(t) & = & \inf\limits_{x, y \in \dom f, \atop \alpha \in (0,1)} \Big\{ {\alpha f(x) + (1-\alpha) f(y) - f(\alpha x + (1-\alpha) y) \over \alpha(1-\alpha)} : \; \| x - y \| \geq t \Big\}.
\ea
\eeq
This modulus satisfies the opposite inequality (see  \cite{VNC})
\beq\label{eq-ModUC}
\ba{rcl}
\mu(\beta t) & \leq & \beta^2 \mu(t), \quad t \in \dom \mu, \; \beta \in [0,1],
\ea
\eeq 
Unfortunately, for non-convex functions this logic does not work and we need to provide inequality (\ref{eq-ModUS}) with a direct proof.

Denote $D = {\rm diam }\, \left(\dom f \right)$, and let $\Gamma = [0,D]$ for $D < +\infty$ and $\Gamma = [0,+\infty)$ otherwise.
In what follows, we assume the following regularity condition.
\BAS\label{ass-Dom}
Function $\hat \mu(\cdot)$ is well defined for all $t \in \Gamma$ and $\lim\limits_{\tau \to +0} {1 \over \tau} \hat \mu(\tau) = 0$.
\EAS

This assumption is not very restrictive. It can be ensured by different standard assumptions on function $f(\cdot)$. This can be seen by the following observation.
\BL\label{lm-Ex}
Denote $\hat \mu_1(t) = \sup\limits_{x, y \in \dom f} \Big\{ \| \nabla f(x) - \nabla f(y)\|_* : \; \| x - y \| \leq t \Big\}$.
Then, for any $t \in \Gamma$, we have
\beq\label{eq-MUp}
\ba{rcl}
\hat \mu(t) & \leq &  \int\limits_0^t \hat \mu_1(\tau) d \tau.
\ea
\eeq
\EL
\proof
Indeed, for all $x, y \in \dom f$ and $\alpha \in (0,1)$, we have
$$
\ba{rcl}
f(x + (1-\alpha)(y-x)) - f(x) & = & (1-\alpha) \int\limits_0^1 \la \nabla f(x + \tau (1-\alpha) (y-x)), y- x \ra d \tau,\\
f(y + \alpha(x-y)) - f(y) & = & \alpha \int\limits_0^1 \la \nabla f(y + \tau \alpha (x-y)), x-y \ra d \tau.
\ea
$$
Hence, 
$$
\ba{rl}
& | \alpha f(x) + (1-\alpha) f(y) - f(\alpha x + (1-\alpha) y) | \\
\\
= &  \alpha (1-\alpha) \Big| \int\limits_0^1 \la \nabla f(x + \tau (1-\alpha) (y-x)) - \nabla f(y + \tau \alpha (x-y)), x-y \ra d \tau \Big| \\
\\
\leq &  \alpha (1-\alpha) t \int\limits_0^1 \| \nabla f(x + \tau (1-\alpha) (y-x)) - \nabla f(y + \tau \alpha (x-y)) \|_* d \tau\\
\\
\leq & \alpha (1-\alpha) t \int\limits_0^1 \hat \mu_1((1-\tau)t) d \tau \; = \; \alpha(1-\alpha) \int\limits_0^t \hat \mu_1(\tau ) d \tau. \QR
\ea
$$

Thus, Assumption \ref{ass-Dom} is ensured, for example, by H\"older condition for the gradient of function $f(\cdot)$, or by uniform boundedness of the norm of its gradients on $\dom f$.

Let us prove the main property of GCB. Our proof is based on the ideas of \cite{VNC}.
\BL\label{lm-HMu}
For any $t \in \Gamma$ and $\beta \in [0,1]$, we have $\hat \mu( \beta t) \geq \beta^2 \hat \mu(t)$.
\EL
\proof
Let us choose first $\half < \beta \leq 1$. For $t \in \Gamma$ and  a small $\epsilon > 0$, let us choose points $u_1, u_2 \in \dom f$ such that $\| u_1 - u_2 \|\leq t$ and
\beq\label{eq-MClose}
\ba{rcl}
| \alpha f(u_1) + (1-\alpha) f(u_2) - f(\alpha u_1 + (1-\alpha) u_2) | & \geq & \alpha (1-\alpha) [ \hat \mu(t) - \epsilon].
\ea
\eeq
Without loss of generality, we can assume that $\alpha \in (0,\half]$.

For $\tau \in [0,1]$, denote $u_{\tau} = \tau u_1 + (1-\tau)u_2$.   Note that $\| u_{\beta} - u_2 \| \leq \beta t$. At the same time, for $\gamma \Def {\alpha \over \beta} < 1$, we have
$$
\ba{rcl}
\gamma u_{\beta} + (1-\gamma) u_2 & = & \gamma [\beta u_1 + (1-\beta) u_2] + (1-\gamma) u_2 \; = \; u_{\alpha}.
\ea
$$
Therefore,
$$
\ba{rcl}
| \beta f(u_1) + (1-\beta) f(u_2) - f(u_{\beta}) | & \leq & \beta(1-\beta) \hat \mu(t),\\
\\
| \gamma f(u_{\beta}) + (1-\gamma) f(u_2) - f(u_{\alpha}) | & \leq & \gamma(1-\gamma) \hat \mu(\beta t).
\ea
$$
Thus,
$$
\ba{rcl}
\alpha (1-\alpha) [ \hat \mu(t) - \epsilon] & \leq & |\alpha f(u_1) + (1-\alpha) f(u_2) - f(u_{\alpha}) |\\
\\
& \leq & \gamma | \beta f(u_1) + (1-\beta) f(u_2) - f(u_{\beta}) |\\
\\
& & 
+ | \gamma f(u_{\beta}) + (1-\gamma) f(u_2) - f(u_{\alpha}) |\\
\\
& \leq & \alpha(1-\beta) \hat \mu(t) + \gamma(1-\gamma) \hat \mu(\beta t).
\ea
$$
Dividing this inequality by $\alpha(\beta - \alpha)$, we get $\beta^{-2} \hat \mu(\beta t) - \hat \mu(t) \geq - {1 - \alpha \over \beta - \alpha} \epsilon$. Choosing $\epsilon \to 0$, we obtain the desired inequality.

Finally, if $\beta \leq \half$, then we can represent it as $\beta = \bar \beta^k$ with some $\bar \beta \in [\half,1)$ and $k \geq 1$, and repeat the above reasoning for $\bar \beta$ recursively $k$ times.
\qed

Let us prove now several useful inequalities for GCB.
\BL\label{lm-Up}
For any $x, y \in \dom f$, we have
\beq\label{eq-FUp}
\ba{rcl}
| f(y) - f(x) - \la \nabla f(x), y - x \ra | & \leq &  \hat \mu(\| y - x \|),
\ea
\eeq
\beq\label{eq-GUp}
\ba{rcl}
|\la \nabla f(x) - \nabla f(y), x - y \ra| & \leq & 2 \hat \mu(\| x - y \|),
\ea
\eeq
\beq\label{eq-FUp2}
\ba{rcl}
| f(y) - f(x) - \la \nabla f(x), y - x \ra | & \leq & 2 \hat \sigma(\| x - y \|), 
\ea
\eeq
where $\hat \sigma(r) \Def \int\limits_0^1 {1 \over \tau} \hat \mu(\tau r) d \tau = \int\limits_0^r {1 \over \tau} \hat \mu(\tau) d \tau$, $r \in \Gamma$.
\EL
\proof
Indeed, in view of definition of $\hat \mu(\cdot)$, we have
$$
\ba{rcl}
\Big| f(y) - f(x) - {1 \over 1 - \alpha} [ f(x + (1-\alpha)(y-x)) - f(x)] \Big| & \leq & \alpha \hat \mu (\| y - x \|).
\ea
$$
Choosing $\alpha \to 1$, we get inequality (\ref{eq-FUp}). Adding two copies of (\ref{eq-FUp}) with $x$ and $y$ interchanged, we get an upper bound for the left-hand side of inequality (\ref{eq-GUp}). Finally, for $r = \| x - y \|$, we have
$$
\ba{rcl}
\Big| f(y) - f(x) - \la \nabla f(x), y - x \ra \Big| & = & \Big| \int\limits_0^1 \la \nabla f(x + \tau (y - x)) - \nabla f(x), y - x \ra d \tau \Big|\\
\\
& = & \Big| \int\limits_0^1 {1 \over \tau} \la \nabla f(x + \tau (y - x)) - \nabla f(x), \tau(y - x) \ra d \tau \Big| \\
\\
& \refLE{eq-GUp} & \int\limits_0^1 {2 \over \tau}  \hat \mu(\tau r) d \tau \; = \; 2 \hat \sigma(r). \QR
\ea
$$

As compared with (\ref{eq-FUp}), inequality (\ref{eq-FUp2}) is weaker.
However, under some natural conditions, function $\hat \sigma(\cdot)$ is proportional to $\hat \mu(\cdot)$. Let us justify its useful properties. 
\BL\label{lm-PropH}
Function $\hat \sigma(\cdot)$ is monotonically increasing on $\Gamma$ and $\lim\limits_{r \to +0} {1 \over r} \hat \sigma(r) = 0$. 
For all $r \in \Gamma$, we have
\beq\label{eq-SigUS}
\ba{rcl}
\hat \sigma(\beta r) & \geq & \beta^2 \hat \sigma(r), \quad \beta \in [0,1],
\ea
\eeq
\beq\label{eq-Comp1}
\ba{rcl}
\hat \mu(r) & \leq & 2 \hat \sigma (r).
\ea
\eeq
Moreover, if $f(\cdot)$ is convex and $\dom f = \E$, then for any $\beta \in [0,1]$, we have
\beq\label{eq-HMuUp}
\ba{rcl}
\hat \mu (\beta r) & \leq & \beta \hat \mu(r). 
\ea
\eeq
Consequently, in this case, function $\hat \sigma(\cdot)$ is convex and for any $r \in \Gamma$ we have
\beq\label{eq-RelHSM}
\ba{rcl}
\hat \sigma(r) & \leq & \hat \mu(r).
\ea
\eeq
\EL
\proof
The first statement of the lemma follows from l'Hopital rule and Assumption \ref{ass-Dom}. If $r \in \Gamma$, then
$\hat \sigma(r) = \int \limits_0^1 {1 \over \tau} \hat \mu(\tau r) d \tau \; \refGE{eq-ModUS} \; \int \limits_0^1 \tau \hat \mu(r) d \tau \; = \; \half \hat \mu(r)$,
and this is (\ref{eq-Comp1}). By the same reason,
$$
\ba{rcl}
\hat \sigma(\beta r) = \int \limits_0^1 {1 \over \tau} \hat \mu(\tau \beta r) d \tau \; \refGE{eq-ModUS} \; \beta^2 \int \limits_0^1 {1 \over \tau} \hat \mu(\tau r) d \tau \; = \; \beta^2 \hat \sigma(r).
\ea
$$

Further, let $f(\cdot)$ be convex and $\dom f = \E$. Let us fix an arbitrarily small $\epsilon > 0$. Then there exist points $u_1, u_2 \in \E$ and $\alpha \in (0,1)$ such that $\| u_1 - u_2 \| \leq t$ and
$$
\ba{rcl}
\alpha f(u_1) + (1-\alpha) f(u_2) - f(u_{\alpha}) & \geq & \alpha(1-\alpha) [\hat \mu(t) - \epsilon],
\ea
$$
where $u_{\alpha} = \alpha u_1 + (1-\alpha)u_2$. Let us choose some $\gamma > 0$ and define
$$
\ba{rcl}
\bar u_i & = & u_i + \gamma(u_i - u_{\alpha}), \quad i = 1, 2.
\ea
$$
Then $\| \bar u_1 - \bar u_1 \| \leq (1+\gamma)t$, and 
$$
\ba{rcl}
\alpha \bar u_1 + (1-\alpha) \bar u_2 & = & \alpha (u_1 + \gamma(u_1 - u_{\alpha}) + (1 - \alpha) (u_2 + \gamma(u_2 - u_{\alpha})) \; = \; u_{\alpha}.
\ea
$$
At the same time, since $f(\cdot)$ is convex, we have 
$f(\bar u_i) \geq f(u_i) + \gamma (f(u_i) - f(u_{\alpha}))$, $i = 1,2$. Therefore,
$$
\ba{rl}
& \hat \mu((1+\gamma)t) \; \geq \; {1 \over \alpha(1-\alpha)} \Big\{ \alpha f(\bar u_1) + (1-\alpha) f(\bar u_2) - f(u_{\alpha}) \Big\}\\
\\
\geq & {1 \over \alpha(1-\alpha)} \Big\{ \alpha [f(u_1) + \gamma (f(u_1) - f(u_{\alpha}))] + (1-\alpha) [ f(u_2) + \gamma (f(u_2) - f(u_{\alpha}))] - f(u_{\alpha}) \Big\}\\
\\
= & {1+\gamma \over \alpha(1-\alpha)} \Big[ \alpha f(u_1) + (1-\alpha) f(u_2) - f(u_{\alpha}) \Big] \; \geq \; (1+\gamma) [ \hat \mu(t) - \epsilon].
\ea
$$
Choosing now $\epsilon>0$ arbitrarily small and denoting $r = (1+\gamma)t$ with $\beta = {1 \over 1+\gamma}$, we get inequality (\ref{eq-HMuUp}). 

We can see now that derivative $\hat \sigma'(t) = {1 \over t} \hat \mu(t)$ is non-decreasing on $\Gamma$:
$$
\ba{rcl}
\hat \sigma '(\beta t) & = & {1 \over \beta t} \hat \mu(\beta t) \; \refLE{eq-HMuUp} \; {1 \over t} \hat \mu(t) \; = \; \hat \sigma'(t), \quad \beta \in (0,1).
\ea
$$
Hence, function $\hat \sigma(\cdot)$ is convex. At the same time, by definition of $\hat \sigma(t)$, we have 
$$
\ba{rcl}
\hat \sigma(r) & = & \int\limits_0^1 {1 \over \tau} \hat \mu(\tau r) d \tau \; \refLE{eq-HMuUp} \; \int\limits_0^1 {1 \over \tau} \cdot \tau  \hat \mu(r) d \tau \; = \; \hat \mu(r). \QF
\ea
$$

As we will see later, the right-hand side of inequality (\ref{eq-FUp2}) is more convenient for the complexity analysis than that of inequality (\ref{eq-FUp}). 

Let us prove now the following important property. Denote $\Gamma_1 = \{ \tau\geq 0: \; \sqrt{\tau} \in \Gamma \}$.
\BL\label{lm-Xi}
On the set $\Gamma_1$, function $\hat \xi(\tau) = \hat \sigma(\sqrt{\tau})$ is concave.
\EL
\proof
Indeed, function $\hat \xi(\cdot)$ is differentiable and $\hat \xi'(\tau) = {1 \over 2 \sqrt{\tau}} \, \hat \sigma'(\sqrt{\tau}) = {1 \over 2 \tau} \hat \mu(\sqrt{\tau})$. Hence, for any $\beta \in (0,1)$, we have
$$
\ba{rcl}
\hat \xi'(\beta \tau) & = & {1 \over 2 \beta \tau} \hat \mu(\sqrt{\beta \tau}) \; \refGE{eq-ModUS} \; {1 \over 2 \tau} \hat \mu(\sqrt{\tau}) \; = \; \hat \xi'(\tau).
\ea
$$
Thus, the derivative $\hat \xi'(\cdot)$ is non-increasing and therefore $\hat \xi(\cdot)$ is concave.
\qed

\BC\label{cor-Xi}
For any values $r, t \in \Gamma$, we have
\beq\label{eq-SigUp}
\ba{rcl}
\hat \sigma(t) & \leq & \hat \sigma(r) - \half \hat \mu(r) + \half \left({t \over r} \right)^2 \hat \mu(r).
\ea
\eeq
\EC
\proof
Let us choose $\tau_r = r^2$ and $\tau = t^2$. Since $\hat \xi(\cdot)$ is concave, we have
$$
\ba{rcl}
\hat \xi(\tau) & \leq & \hat \xi(\tau_r) + \hat \xi'(\tau_r)(\tau - \tau_r) \; = \; \hat \sigma(r) + {1 \over 2 r^2} \hat \mu (r) (t^2 - r^2),
\ea
$$
and this is (\ref{eq-SigUp}). 
\qed

For $x, y \in \dom f$, we introduce the Bregmann distance between $x$ and $y$ as usual:\footnote{For nonconvex $f(\cdot)$, this "distance" is sign-indefinite.}
$$
\ba{rcl}
\beta_f(x,y) & = & f(y) - f(x) - \la \nabla f(x), y - x \ra.
\ea
$$
By inequalities (\ref{eq-FUp2}), we come to the following statement.
\BT\label{th-BBounds}
For any $r \in \Gamma$ and $x, y \in \dom f$, we have the following bounds:
\beq\label{eq-BBound}
\ba{rcl}
| \beta_f(x,y) | & \leq & \delta^+_f(r) + \half L_f(r) \| y - x \|^2,
\ea
\eeq
where
$$
\ba{rcl}
\delta^+_f(r) & \Def & 2 \hat \sigma(r) - \hat \mu(r) \refGE{eq-Comp1} \; 0, \quad L_f(r) \; \Def \; {2 \over r^2} \hat \mu(r).
\ea
$$
\ET

\BE\label{ex-Quad}
For $A = A^T \in \R^{n \times n}$, define $f(x) = \half \la A x, x \ra$. Then for all $x , y \in \E$ and $\alpha \in (0,1)$, we have
$$
\ba{rcl}
\alpha f(x) + (1- \alpha) f(y) - f(\alpha x + (1-\alpha)y) & = & \half \alpha(1-\alpha) \la A(x-y),x - y \ra.
\ea
$$
Therefore, $\hat \mu(r) = \half L r^2$ with $L = \max\left\{ \lambda_{\max}(A), - \lambda_{\min}(A) \right\}$. Thus,
$$
\ba{rcl}
\hat \sigma(r) & = & \int\limits_0^r {1 \over \tau} \hat \mu(\tau) d \tau \; = \; {1 \over 4} L r^2, \quad \delta^+_f(r) \; = \; 0, \quad L_f(r) \; \equiv \; L.
\ea
$$
Hence, in our example, inequality (\ref{eq-BBound}) results in the following exact bound
$$
\ba{rcl}
|f(y) - f(x) - \la \nabla f(x), y - x \ra | & \leq & \half L \| y - x \|^2. \QF
\ea
$$
\EE

\section{Complexity of universal gradient methods}\label{sc-CompGM}
\SetEQ

In order to use the bound (\ref{eq-BBound}) for analysing the universal gradient methods, we need to justify some regularity of the functions $\delta^{+}_f(r)$ and $L_f(r)$. For that, we need the following simple lemma.
\BL\label{lm-Grow}
Let function $\psi(\cdot)$ be defined on $\Gamma$. Assume that for some $r_0 \in \Gamma$ there exist constants $a, b > 0$ such that for any $r \in \Gamma$, $r \geq r_0$, and $\alpha \in [0,b]$ with $r+\alpha \in \Gamma$, we have $\psi(r+\alpha) \geq \psi(r) - a \alpha^2$. Then $\psi(r) \geq \psi(r_0)$ for all $r \geq r_0$, $r \in \Gamma$.
\EL
\proof
Let us fix some $\bar r > r_0$, $\bar r \in \Gamma$. Let us choose $N \geq 1$ big enough to have $\alpha = {\bar r - r_0 \over N} \leq b$. Consider the points $r_k = r_0 + k \alpha$. Then $\psi(r_{k+1}) \geq \psi(r_k) - a \alpha^2$, $k = 0, \dots, N-1$. Hence,
$\psi(\bar r) = \psi(r_N) \geq \psi(r_0) - N a \alpha^2 = \psi(r_0) - {a \over N} (\bar r - r_0)^2$. Thus, choosing $N \to \infty$, we get the result.
\qed

Now we can prove the following statement.
\BL\label{lm-CGrow}
On the set $\Gamma$, function $\delta^{+}_f(\cdot)$ is non-decreasing, and function $L_f(\cdot)$ is non-increasing.
\EL
\proof
Consider $r \in \Gamma$ and $\alpha > 0$ such that $r+\alpha \in \Gamma$. Then
$$
\ba{rcl}
\delta^+_f(r+\alpha) - \delta^+_f(r) & = & 2 \int\limits_r^{r+\alpha} {1 \over \xi} \hat \mu(\xi) d \xi - \Big[ \hat \mu(r+\alpha) - \hat \mu(r) \Big]\\
\\
& \refGE{eq-ModUS} & {2\alpha \hat \mu(r) \over r + \alpha} -  \Big[ \left(1 + {\alpha \over r} \right)^2 - 1\Big]\hat \mu(r) \; = \; - \alpha^2 \Big[ { 2 \over r(r+\alpha)} + {1 \over r^2} \Big] \hat \mu(r).
\ea
$$
Hence, in view of Lemma \ref{lm-Grow}, function $\delta^+_f(\cdot)$ is non-decreasing. 

Finally, for $ r \in \Gamma$ and $\beta \in (0,1)$, we have
$$
\ba{rcl}
L_f(\beta r) = {2 \over \beta^2 r^2} \hat \mu(\beta r) \; \refGE{eq-ModUS} \; L_f(r).
\ea
$$
Hence, function $L_f(\cdot)$ is non-increasing on $\Gamma$.
\qed

In the main part of this paper, we present a performance analysis of universal variants of first-order methods, which were developed in \cite{UGM} for solving the
following optimization problem:\footnote{For integrity of presentation, in Section \ref{sc-NonConv} we analyze a simple Universal Gradient Method for nonconvex problems.}
\beq\label{prob-Comp}
\min\limits_{x \in \dom \Psi} \Big[ \tilde f(x) = f(x) + \Psi(x) \Big],
\eeq
where $f(\cdot)$ is a convex function and $\Psi(\cdot)$ is a simple closed convex function with $\dom \Psi \subseteq \dom f \subseteq \E$, $\inter (\dom \Psi) \neq \emptyset$. We assume that at any $x \in \dom \Psi$ there exists at least one subgradient $\nabla f(x) \in \partial f(x)$. If $f(\cdot)$ is differentiable at $x$, then this is the usual gradient.

Let us explain our approach to the performance analysis of universal gradient methods under
simplifying assumption that $f(\cdot)$ is convex and $\dom f = \dom \Psi = \E$. Then, in view of Lemma \ref{lm-PropH}, we have 
\beq\label{eq-DelMuH}
\ba{rcl}
\delta_f^+(r) \; \refLE{eq-RelHSM} \; \hat \mu(r), \quad r \in \Gamma.
\ea
\eeq

For such methods, we usually get an estimate for the rate of convergence in terms of $\delta^+_f(r)$ and $L_f(r)$. Assume, for example, that for a hypothetical method ${\cal M}$, we can justify the following rate of convergence:
\beq\label{eq-fbound}
\ba{rcl}
\tilde f(x_k) - \tilde f(x^*) & \leq & \delta^+_f(r) + {1 \over 2 k} L_f(r) r_0^2 \; \refLE{eq-DelMuH} \; \hat \mu(r) + {1 \over k} \hat \mu(r) \left({r_0 \over r}\right)^2,
\ea
\eeq
where $k \geq 1$ is the number of iterations of the method, $r \in \Gamma$, and $r_0$ is an upper bound for the distance between the starting point and an arbitrary optimal solution $x^*$. 

For universal gradient methods, the only input parameter is the desired accuracy $\epsilon>0$ of the solution of problem (\ref{prob-Comp}) in function value. Let us choose the radius $r_{\epsilon} \in \Gamma$ from the condition
$$
\ba{rcl}
\hat \mu(r_{\epsilon}) & = & \half \epsilon.
\ea
$$
For further reasoning, we need the following definition.
\BD\label{def-CompGauge}
\underline{\em Complexity gauge} $s_f(\cdot)$ for function $f(\cdot)$ is defined by the following equation:
\beq\label{def-CGauge}
\ba{rcl}
\hat \mu_f(s_f(\epsilon)) & = & \epsilon, \quad \epsilon \geq 0.
\ea
\eeq
In other words, the complexity gauge is an inverse function of the Global Curvature Bound of function $f(\cdot)$.
\qed
\ED

Then, the number of iterations $k$, which is sufficient for getting an $\epsilon$-solution of the problem, satisfies the following inequality:
\beq\label{eq-UComp}
\ba{rcl}
k & \geq & {2 \over \epsilon} \hat \mu(r_{\epsilon}) \left({r_0 \over r_{\epsilon}}\right)^2 \; = \; \left({r_0 \over r_{\epsilon}}\right)^2 \; = \; {r^2_0 \over s^2_f\left(\epsilon\over 2\right)}.
\ea
\eeq
We call this complexity bound {\em inverse} since it is written in terms of inverse function $s_f(\cdot)$.

For some universal methods, it is more convenient to rewrite the bound (\ref{eq-UComp}) in its {\em direct} form, which is similar to the usual rate of convergence. Indeed, inequality (\ref{eq-UComp}) tells us that any number of iterations $k$ satisfying the inequality 
\beq\label{eq-UComp1}
\ba{rcl}
2 \hat \mu\left( {r_0 \over \sqrt{k}}\right) & \leq & \epsilon
\ea
\eeq
is sufficient for getting $\epsilon$-solution of our problem. Inequalities of this type we call the {\em direct complexity bounds}.

Thus, the Global Curvature Bound and the complexity gauge provide us with universal complexity bounds for solving the problem (\ref{prob-Comp}). Note that these objects describe the interaction of the epigraph of function $f(\cdot)$ with the norm $\| \cdot \|$. They do not depend on artificially introduced problem classes containing $f(\cdot)$. Moreover, the bounds~(\ref{eq-UComp}) and~(\ref{eq-UComp1}) are sharper than any bound obtained from a parametric problem class related to the level of smoothness of the objective function.

\BE\label{ex-CN}
Let $f \in \C^{1,\nu}_{L_{\nu}}(\E)$. Then
$\hat \mu(t) \refLE{eq-MUp} L_{\nu} t^{1+\nu}$. Hence
$\epsilon = \hat \mu(s_f(\epsilon)) \; \leq \; L_{\nu} s^{1+\nu}_f(\epsilon)$,
and we conclude that 
$$
\ba{rcl}
s_f(\epsilon) & \geq & \left({\epsilon \over L_{\nu}} \right)^{1 \over 1+ \nu}.
\ea
$$
Hence, for the complexity bound (\ref{eq-UComp}), we have
\beq\label{eq-SBound}
\ba{rcl}
{r^2_0 \over s^2_f\left(\epsilon\over 2\right)} \leq r^2_0 \left({2L_{\nu} \over \epsilon}\right)^{2 \over 1+ \nu}.
\ea
\eeq
Note that the right-hand side of this inequality corresponds to the rate of convergence
$$
\ba{rcl}
f(x_k) - f(x^*) & \leq & {2 L_{\nu} r_0^{1+\nu} \over k^{1+\nu \over 2}}, \quad k \geq 1,
\ea
$$
which is typical for the variants of Gradient Method. In view of (\ref{eq-SBound}), the bound (\ref{eq-UComp}) is not worse than the complexity bound for any class $\C^{1,\nu}_{L_{\nu}}(\E)$ containing $f(\cdot)$. \qed
\EE

Let us show how we can treat some non-standard problem classes.
\BE\label{ex-Sum}
Let $f_0 \in \C^{1,0}(\E)$ and $f_1 \in \C^{1,1}(\E)$. For problem (\ref{prob-Comp}), consider the objective function 
\beq\label{eq-FAdd}
\ba{rcl}
f(x) & = & f_0(x) + f_1(x).
\ea
\eeq 
Since $\hat \mu_{f_0+f_1}(r) \leq \hat \mu_{f_0}(r) + \hat \mu_{f_1}(r)$, $r \geq 0$, we have
\beq\label{eq-Sum}
\ba{rcl}
\hat \mu_f(r) & \leq & \hat \mu_{f_0}(r) + \hat \mu_{f_1}(r) \; \refLE{eq-MUp} \; L_0 r + \half L_1 r^2.
\ea
\eeq
Therefore, $s_f(\epsilon) \geq r_{\epsilon}$, where $r_{\epsilon}$ is the positive root of the following equation:
$$
\ba{rcl}
L_0 r + \half L_1 r^2 & = & \epsilon.
\ea
$$
This is $r_{\epsilon} = {2 \epsilon \over L_0 + \sqrt{2 \epsilon L_1 + L_0^2}}$. Thus, a sufficient number of steps for solving problem (\ref{prob-Comp}) with $f(\cdot)$ given by (\ref{eq-FAdd}) up to accuracy $\epsilon>0$, satisfies the following bound:
\beq\label{eq-KSuff}
\ba{rcl}
k & \geq & r_0^2 \Big[ {L_0 \over \epsilon} + \sqrt{ {L_1 \over \epsilon} + {L_0^2 \over \epsilon^2}} \Big]^2 \; = \; O\left(r_0^2 \Big[ {L_0^2 \over \epsilon^2} + {L_1 \over \epsilon} \Big]\right).
\ea
\eeq
(In accordance to (\ref{eq-UComp}), our conclusion is based on $r_{\epsilon/2}$.) 
\qed
\EE

Thus, for analyzing the method ${\cal M}$ on a new problem class, it is not necessary to perform a new investigation of the method itself. We just need to derive a new lower bound for the complexity gauge $s_f(\cdot)$.

\section{Universal Gradient Method for Nonconvex \\ Problems}\label{sc-NonConv}
\SetEQ

We start our presentation of Universal Gradient Methods from a simple method for solving the following constrained minimization problem:
\beq\label{prob-Const}
f_* \; \Def \; \min\limits_{x \in Q} \; f(x) \; > \; - \infty,
\eeq
where $Q \subseteq \dom f \subseteq \E$ is a {\em simple} closed convex set, and a general objective function $f(\cdot)$ satisfies Assumption \ref{ass-Dom}. Let us measure distances in $\E$ by the Euclidean norm
$$
\ba{rcl}
\| x \| & = & \la B x, x \ra^{1/2}, \quad x \in \E,
\ea
$$
where $B = B^*: \; \E \to \E^*$ is positive definite. We assume that the objective function in problem (\ref{prob-Const}) is below bounded by some value $f_* > - \infty$.

For approaching the saddle points of problem (\ref{prob-Const}), let us define the {\em Gradient Mapping} (e.g. Section 2.2.4 in \cite{Lect}):
\beq\label{def-GM}
\ba{rcl}
{\cal T}_M(\bar x) & = & \arg\min\limits_{x \in Q} \Big\{ f(\bar x) + \la \nabla f(\bar x), x - \bar x \ra + \half M \| x - \bar x \|^2 \Big\},\\
\\
g_M(\bar x) & = & M( \bar x - {\cal T}_M(\bar x)),
\ea
\eeq
where $\bar x \in Q$ and $M > 0$. The point $T = {\cal T}_M(\bar x)$ satisfies the following optimality condition:
\beq\label{eq-OptT}
\ba{rcl}
\la \nabla f(\bar x) + M B(T-\bar x), x - T \ra & \geq & 0, \quad \forall x \in Q.
\ea
\eeq
Applying this inequality to $x = \bar x$, we get 
\beq\label{eq-GMBound}
\ba{rcl}
\la \nabla f(\bar x), \bar x - T \ra & \geq & {1 \over M} \| g_M(\bar x) \|_*^2.
\ea
\eeq
Thus, $\| g_M(\bar x) \|* \leq \| \nabla f(\bar x)||_*$ for any $M >0$. In the case $Q \equiv \E$, we have $g_M(\bar x) = \nabla f(\bar x)$.

For non-convex functions, it is often better to use $\hat \sigma(\cdot)$ instead of $\hat \mu(\cdot)$. Let us define for it the corresponding inverse function $\hat s_f(\cdot)$ as follows:
\beq\label{def-HatS}
\ba{rcl}
\hat \sigma_f(\hat s_f(\epsilon)) & = & \epsilon, \quad \epsilon \in \hat \sigma(\Gamma).
\ea
\eeq
Let us prove the following useful inequality.
\BL\label{lm-HSig}
For any $t \in \hat \sigma(\Gamma)$ and $\beta \in [0,1]$ we have
\beq\label{eq-HSig}
\ba{rcl}
\hat s_f(\beta t)& \leq & \sqrt{\beta} \hat s_t(t).
\ea
\eeq
\EL
\proof
Let $r = \hat s_f(t)$. Then $\hat \sigma(r) = t$. Therefore, by (\ref{eq-SigUS}), we have
$$
\ba{rcl}
\beta t & = & \beta \hat \sigma(r) \leq \hat \sigma\left(\sqrt{\beta} r\right).
\ea
$$
Hence, $\hat s_f(\beta t) \leq \sqrt{\beta} r = \sqrt{\beta} \hat s_f(t)$.
\qed

Further, for any $t \in 2 \hat \sigma(\Gamma)$, define 
\beq\label{def-HGammaF}
\ba{c}
\mbox{\fbox{$\hat \gamma_f(t) = {2 t \over \hat s_f^2\left({t \over 2}\right)}$}}
\ea
\eeq
Note that this function is non-increasing:
\beq\label{eq-HGNon} 
\ba{rcl}
\hat \gamma_f(\beta t) & = & {2 \beta t \over \hat s_f^2 \left({\beta t \over 2} \right)} \; \; \refGE{eq-HSig} \; {2 t \over \hat s_f^2 \left({t \over 2} \right)} \; = \; \hat \gamma_f(t), \quad t \in 2 \hat \sigma(\Gamma), \; \beta \in (0,1].
\ea
\eeq

We can prove now a convenient upper bound for the Bregmann Distance. 
\BL\label{lm-ModS}
For any $\epsilon \in 2 \hat \sigma(\Gamma)$ and any $M \geq \hat \gamma_f(\epsilon)$, we have
\beq\label{eq-ModSUp}
\ba{rcl}
|\beta_f(x,y)| & \leq & \half M \| y - x \|^2 + {\epsilon \over 2}, \quad x, y \in \dom f.
\ea
\eeq
\EL
\proof
Indeed, for $r=\hat s_f\left({\epsilon \over 2} \right)$, we have $\delta^+_f(r) \leq \hat \sigma(r) = {\epsilon \over 2}$ and
$$
\ba{rcl}
L_f(r) & = & {2 \over r^2} \hat \mu(r) \; \refLE{eq-Comp1} \; {4 \over r^2} \hat \sigma(r)
\; = \; {2 \epsilon \over \hat s_f^2 \left({\epsilon \over 2}\right)} \; = \; \hat \gamma_f(\epsilon).
\ea
$$
Hence, (\ref{eq-ModSUp}) follows from (\ref{eq-BBound}).
\qed

Lemma \ref{lm-ModS} allows us to estimate the tightness of the first-order optimality condition at point $T = {\cal T}_M(\bar x)$.
\BT\label{th-Crit}
Let $\bar x \in Q$, $\epsilon \in 2 \sigma(\Gamma)$, $M \geq \hat \gamma_f(\epsilon)$, and $T = {\cal T}_M(\bar x)$. Then 
\beq\label{eq-Crit}
\ba{rcl}
\la \nabla f(T) , x - T \ra & \geq & - {3 \over 2} \epsilon - {2 \over M} \| g_M(\bar x) \|_*^2 - 2M \| x - T \|^2, \quad x \in Q.
\ea
\eeq
\ET
\proof
Denote $r = \| T - \bar x \|$. Applying inequality (\ref{eq-ModSUp}) three times, for any $x \in Q$, we have
$$
\ba{rl}
& \la \nabla f(T) , x - T \ra \; \geq \; f(x) - f(T) - {M \over 2} \| x - T \|^2 - {\epsilon \over 2}\\
\\
\geq & f(\bar x) + \la \nabla f(\bar x), x - \bar x \ra - f(T) - {M \over 2} \| x - \bar x \|^2  - {M \over 2} \| x - T \|^2 - \epsilon\\
\\
\refGE{eq-OptT} &
f(\bar x) + \la \nabla f(\bar x), T - \bar x \ra + \la g_M(\bar x), x - T \ra - f(T) - {M \over 2} \| x - \bar x \|^2  - {M \over 2} \| x - T \|^2 - \epsilon\\
\\
\geq & \la g_M(\bar x), x - T \ra - {M \over 2} r^2 - {M \over 2} \| x - \bar x \|^2  - {M \over 2} \| x - T \|^2 - {3 \over 2}\epsilon.
\ea
$$
It remains to note that
$$
\ba{rl}
& \la g_M(\bar x), T-x \ra + {M \over 2} \left(r^2 + \| x - \bar x \|^2  + \| x - T \|^2\right)\\
\\
\leq & \la g_M(\bar x), T-x \ra + {3M \over 2} \left(r^2 + \| x - T \|^2\right) \; \leq \; 2 M \left(r^2 + \| x - T \|^2 \right). \QF
\ea
$$

Consider the following optimization scheme.
\beq\label{met-GGM}
\ba{|l|}
\hline \\
\hspace{10ex} \mbox{\bf Universal General Gradient Method (GGM)}\\
\\
\hline \\
\mbox{{\bf Initialization.} Choose a point $x_0 \in Q$, some $M_0 > 0$, and accuracy $\delta > 0$.}\\
\\
\mbox{{\bf For $k \geq 0$ do:}}\\
\\
\ba{rl}
1. & \mbox{Starting from $i = 0$, perform the following operations:}\\
&\ba{rl} \\
\mbox{a)} & \mbox{Compute $x_k^i = {\cal T}_{2^{i}M_k}(x_k)$. Define $g^i_k = 2^{i}M_k(x_k - x_k^i)$.}\\
\\
\mbox{b)}& \mbox{If $\| g_k^i \|_* \leq \delta$, then \underline{\sc Stop}. Define the output $\bar x = x_k^i$, $\bar M = 2^i M_k$.}\\
\\
\mbox{c)} & \mbox{Check validity of inequality $f(x_k)- f(x^+_k) \geq 2^{i-2}M_k \| x_k^i - x_k \|^2$.}\\
\\
\mbox{d)} & \mbox{If yes, then set $i_k = i$ and go to Step 2. Otherwise, set $i:=i+1$.}
\ea\\
\\
2. & \mbox{Set $x_{k+1} = x^{i_k}_k$ and $M_{k+1} = 2^{i_k-1} M_k$.} \\
\ea\\
\\
\hline
\ea
\eeq

Denote $\Delta_0 = f(x_0)-f_*$.
\BT\label{th-GGM}
Let $M_0 \leq \hat \gamma_f(\Delta_0)$ and $N \geq 1$ satisfies the following inequality:
\beq\label{eq-GGMStop}
\ba{rcl}
\Delta_0 & > &  N \hat \sigma_f \left( 2^{5/2} {\Delta_0 \over N \delta} \right).
\ea
\eeq
Then the termination criterion of method (\ref{met-GGM}) will be activated at some iteration $\bar k \leq N$. For all $k$, $0 \leq k \leq \bar k$, we have
\beq\label{eq-UpMG}
\ba{rcl}
M_k & \leq & \hat \gamma_f(\epsilon_N),
\ea
\eeq
where $\epsilon_N ={2\Delta_0 \over N}$.
\ET
\proof
Note that the inequality (\ref{eq-GGMStop}) can be rewritten as $\delta \hat s_f\left({\epsilon_N \over 2}\right) > 2^{3/2} \epsilon_N$. In other words, we have
\beq\label{eq-ED}
\ba{rcl}
{\delta^2 \over \epsilon_N} & > & {8 \epsilon_N \over \hat s_f^2\left({\epsilon_N \over 2}\right)} \; \refEQ{def-HGammaF} \; 4 \hat \gamma_f(\epsilon_N).
\ea
\eeq

First of all, we need to prove that method (\ref{met-GGM}) is well defined for all $k = 0, \dots, \bar k$. Let us prove by induction that
\beq\label{eq-UpM}
\ba{rcl}
M_k & \leq & \hat \gamma_f(\epsilon_N), \quad 0 \leq k \leq \bar k.
\ea
\eeq
Indeed, for $k=0$ we have $M_0 \leq \hat \gamma_f(\Delta_0) \refLE{eq-GNon} \hat \gamma_f(\epsilon_N)$. Assume now that (\ref{eq-UpM}) is true for some $k \geq 0$. If $i_k = 0$, then $M_{k+1} = \half M_k \leq \half \hat \gamma_f(\epsilon_N)$. 

Let us assume that $i_k > 0$. Then there are two possibilities. Firstly, the inequality of Step 1c) can be satisfied for some $i > 0$ such that $2^i M_k \leq \hat \gamma_f(\epsilon_N)$. In this case, $M_{k+1} \refLE{eq-UpM} \half \hat \gamma_f(\epsilon_N)$. Or, for some $i\geq 1$ we will have eventually $2^{i-1} M_k \leq \hat \gamma_f(\epsilon_N) \leq 2^i M_k$. Then
$$
\ba{rcl}
f(x_k) - f(x_k^i) & \refGE{eq-ModSUp} & \la \nabla f(x_k), x_k - x_k^i \ra - {2^iM_k \over 2} \| x_k - x_k^i \|^2 - \half \epsilon_N\\
\\
& \refGE{eq-OptT} & {2^iM_k \over 2} \| x_k - x_k^i \|^2 - \half \epsilon_N.
\ea
$$
If at this moment the stopping criterion of Step 1b) does not work, then
$$
\ba{rcl}
2^i M_k \| x_k - x^i_k \|^2 & = & {1 \over 2^i M_k} \| g_k^i \|^2 \; \geq \; {\delta^2 \over 2 \hat \gamma_f(\epsilon_N)} \; \refGE{eq-ED}  \; 2 \epsilon_N.
\ea
$$
This implies that $f(x_k) - f(x_k^i) \geq 2^{i-2} M_k \| x_k - x_k^i \|^2$ and the condition of Step 1c) is satisfied. As we have seen, then $M_{k+1} = 2^{i-1} M_k \leq \hat \gamma_f(\epsilon_N)$. Thus, the relation (\ref{eq-UpM}) is proved.

Let us assume now that the stopping criterion of Step 1b) was not activated during all iterations $0 \leq k \leq N$. Then, in all cases
$$
\ba{rcl}
f(x_k) - f(x_{k+1}) & \geq & {1 \over 8 M_{k+1}} \| g_k^{i_k} \|^2_*  \; \geq \; {\delta^2 \over 8 \hat \gamma_f(\epsilon_N)} \; \refGT{eq-ED} \; \half \epsilon_N \; = \; {\Delta_0 \over N}.
\ea
$$
Summing up these inequalities for $k = 0, \dots, N-1$, we get
$$
\ba{rcl}
\Delta_0 & \geq & f(x_0) - f(x_N) \; > \; \Delta_0.
\ea
$$
This is a contradiction, which proves the last statement of the theorem.
\qed

\BE\label{ex-Holder}
Let us show how the bound (\ref{eq-GGMStop}) works for functions with H\"older continuous gradients:
\beq\label{def-Hold}
\ba{rcl}
\| \nabla f(x) - \nabla f(y) \|_* & \leq & L_{\nu} \| x - y \|^{\nu}, \quad \forall x, y \in Q,
\ea
\eeq
where $\nu \in (0,1]$ is the H\"older parameter and $L_{\nu} > 0$. In accordance to Lemma \ref{lm-Ex}, we have 
$$
\ba{rcl}
\hat \mu(t) & \leq & {L_{\nu} \over 1+ \nu} t^{1+\nu},
\ea
$$
and $\hat \sigma(t) \leq  {L_{\nu} \over (1+ \nu)^2} t^{1+\nu}$. Hence, the sufficient condition for inequality (\ref{eq-GGMStop}) is as follows:
$$
\ba{rcl}
\Delta_0 & > & {N L_{\nu} \over (1+\nu)^2} \left( 2^{5/2} {\Delta_0 \over N \delta} \right)^{1+\nu}.
\ea
$$
This is $N > \left[ {L_{\nu} \over (1+\nu)^2} \left( {2^{5/2} \over \delta} \right)^{1+\nu} \right]^{1/\nu} (f(x_0) - f_*)$.
\qed
\EE

\section{Simple universal gradient methods}\label{sc-UMSimple}
\SetEQ

In this section, we estimate performance of universal gradient methods from \cite{UGM} as applied to problem (\ref{prob-Comp}). We will write our bounds in terms of complexity gauge of function $f(\cdot)$. For simplicity, we assume that $\dom f = \E$. Under this assumption, we can prove the following useful bounds.
\BL\label{lm-GBound}
Let $\dom f = \E$. Then for any $t \geq 0$ and $\beta \in [0,1]$, we have
\beq\label{eq-GBound}
\ba{rcl}
\beta s_f(t) & \leq & s_f(\beta t) \; \leq \; \sqrt{\beta} s_f(t).
\ea
\eeq
\EL
\proof
Let $r = s_f(t)$. Then $\hat \mu(r) = t$. Therefore, we have $\beta t = \beta \hat \mu(r) \refGE{eq-HMuUp} \hat \mu(\beta r)$. This means, that
$s_f(\beta t) \geq \beta r = \beta s_f(t)$. For the second, inequality, we use Lemma \ref{lm-HMu}:
$$
\ba{rcl}
\beta t & = & \beta \hat \mu(r) \leq \hat \mu\left(\sqrt{\beta} r\right).
\ea
$$
Hence, $s_f(\beta t) \leq \sqrt{\beta} r = \sqrt{\beta} s_f(t)$.
\qed

For convex set $\dom \Psi$, we need to
introduce a prox-function $d(\cdot)$ . It is a differentiable convex function, which satisfies the following condition:
\beq\label{def-SConv}
\ba{rcl}
d(y) & \geq & d(x) + \la \nabla d(x), y - x \ra + \half \| y - x \|^2, \quad x \in \inter (\dom \Psi), \; y \in \dom \Psi.
\ea
\eeq
Without loss of generality, we assume that for $x_0 \Def \argmin\limits_{x \in \dom \Psi} d(x)$ we have $d(x_0) = 0$. Thus, 
\beq\label{eq-DGrow}
\ba{rcl}
d(x) & \geq & \half \| x - x_0 \|^2, \quad x \in \dom \Psi.
\ea
\eeq
Clearly, $\beta_d(x,x) = 0$ and
\beq\label{eq-BGrow}
\ba{rcl}
\beta_d(x,y) & \geq & \half \| x - y \|^2, \quad x \in \inter (\dom \Psi), \; y \in \dom \Psi.
\ea
\eeq
Now we can define the {\em Bregman mapping} as follows:
\beq\label{def-BMap}
\ba{rcl}
{\cal B}_M(x) & = & \argmin\limits_{y \in \dom \Psi} \Big\{ \psi_M(x,y) \Def f(x) + \la \nabla f(x), y - x \ra + M \beta_d(x,y) + \Psi(y) \Big\}.
\ea
\eeq
Denote $\psi^*_M(x) = \psi_M(x, {\cal B}_M(x))$. In what follows, we assume that the Bregman mapping can be computed in a closed form. The first-order optimality condition for the point $T = {\cal B}_M(x)$ is as follows:\footnote{As compared with \cite{UGM}, we use another form of the first-order optimality condition, which does not assume subdifferentiability of $\Psi(\cdot)$ at point $T$.}
\beq\label{eq-FOP}
\ba{rcl}
\la \nabla f(x) + M ( \nabla d(T) - \nabla d(x)), y - T \ra + \Psi(y) & \geq & \Psi(T), \quad y \in \dom \Psi.
\ea
\eeq

Let us present now an analogue of Lemma 2 in \cite{UGM}. 
For any $t > 0$, define 
\beq\label{def-GammaF}
\ba{c}
\mbox{\fbox{$\gamma_f(t) = {t \over s_f^2\left({t \over 2}\right)}$}}
\ea
\eeq
Note that this function is non-increasing:
\beq\label{eq-GNon} 
\ba{rcl}
\gamma_f(\beta t) & = & {\beta t \over s_f^2 \left({\beta t \over 2} \right)} \; \; \refGE{eq-GBound} \; {t \over s_f^2 \left({t \over 2} \right)} \; = \; \gamma_f(t), \quad t >0, \; \beta \in (0,1].
\ea
\eeq
\BL\label{lm-MBig}
For any $\epsilon > 0$ and any $M \geq \gamma_f(\epsilon)$, we have
\beq\label{eq-FEUp}
\ba{rcl}
f(y) & \leq & f(x) + \la \nabla f(x), y - x \ra + \half M \| y - x \|^2 + {\epsilon \over 2}, \quad x, y \in \dom \Psi.
\ea
\eeq
Therefore,
\beq\label{eq-FPlus}
\ba{rcl}
\tilde f({\cal B}_M(x)) & \leq & \psi^*_M(x) + {\epsilon \over 2}.
\ea
\eeq
\EL
\proof
Indeed, for $r=s_f\left({\epsilon \over 2} \right)$, we have $\delta^+_f(r) \refLE{eq-DelMuH} \hat \mu(r) = {\epsilon \over 2}$ and
$$
\ba{rcl}
L_f(r) & = & {2 \over r^2} \hat \mu(r) \; = \; {\epsilon \over s_f^2 \left({\epsilon \over 2}\right)^2} \; = \; \gamma_f(\epsilon).
\ea
$$
Hence, (\ref{eq-FEUp}) follows from (\ref{eq-BBound}). Using this inequality, for $x_+ = {\cal B}_M(x)$, we have
$$
\ba{rcl}
f(x_+) & \leq & {\epsilon \over 2} + f(x) + \la f(x), x_+ - x \ra + \half M \| x_+ - x \|^2 \\
\\
& \refLE{eq-DGrow} & {\epsilon \over 2} + f(x) + \la f(x), x_+ - x \ra + M d(x,x_+).
\ea
$$
Therefore, $\tilde f(x_+) = f(x_+) + \Psi(x_+) \leq \psi^*_M(x) + {\epsilon \over 2}$.
\qed

Consider the following optimization scheme \cite{UGM}.
\beq\label{met-PGM}
\ba{|l|}
\hline \\
\hspace{10ex} \mbox{\bf Universal Primal Gradient Method (PGM)}\\
\\
\hline \\
\mbox{{\bf Initialization.} Choose $L_0 > 0$ and accuracy $\epsilon > 0$.}\\
\\
\mbox{{\bf For $k \geq 0$ do:}}\\
\\
\ba{rl}
1. & \mbox{Find the smallest $i_k \geq 0$ such that for $x_k^+ \Def {\cal B}_{2^{i_k}L_k}(x_k)$ we have}\\
\\
& f(x_k^+) \leq f(x_k) + \la \nabla f(x_k), x_k^+ - x_k \ra + 2^{i_k-1}L_k \| x_k^+ - x_k \|^2 + {\epsilon \over 2}.\\
\\
2. & \mbox{Set $x_{k+1} = {\cal B}_{2^{i_k}L_k}(x_k)$ and $L_{k+1} = 2^{i_k-1} L_k$.}\\
\ea\\
\\
\hline
\ea
\eeq
Denote $S_k = \sum\limits_{i=0}^{k} {1 \over L_{i+1}}$ and $\tilde f_k = {1 \over S_k} \sum\limits_{i=0}^{k} {1 \over L_{i+1}} \tilde f(x_i)$.
\BT\label{th-PGM}
Let $f(\cdot)$ satisfy Assumption \ref{ass-Dom} and $L_0 \leq \gamma_f(\epsilon)$. Then, for all $k \geq 0$, we have $L_k \leq \gamma_f(\epsilon)$. Moreover, for all $y \in \dom \Psi$, it holds
\beq\label{eq-RatePGM}
\ba{rcl}
\tilde f^*_k & \leq & {1 \over S_k} \sum\limits_{i=0}^k {1 \over L_{i+1}} [ f(x_i) + \la \nabla f(x_i), y - x_i \ra] + \Psi(y) + {\epsilon \over 2} + {1 \over S_k} \beta(x_0,y).
\ea
\eeq
Therefore, $\tilde f^*_k - \tilde f(x^*) \leq {\epsilon \over 2} + {2 \over k+1} \gamma_f(\epsilon) \beta(x_0,x^*)$.
\ET
\proof
Our proof is very similar to the proof of Theorem 1 in \cite{UGM}. However, since we use another first-order optimality condition (\ref{eq-FOP}), we present it in full details.

In view of inequality (\ref{eq-BBound}), the line-search procedure of Step 1 in method (\ref{met-PGM}) is well defined and
\beq\label{eq-Bound2L}
\ba{rcl}
2 L_{k+1} & = & 2^{i_k} L_k \; \leq \; 2 \gamma_f(\epsilon).
\ea
\eeq

Let us fix an arbitrary point $y \in \dom \Psi$. Denote $r_k \Def \beta(k_k,y)$. Then
$$
\ba{rcl}
r_{k+1}(y) & = & d(y) - d(x_{k+1}) - \la \nabla d(x_{k+1}), y - x_{k+1} \ra\\
\\
& \refLE{eq-FOP} & d(y) - d(x_{k+1}) - \la \nabla d(x_k), y - x_{k+1} \ra \\
\\
& & + {1 \over 2 L_{k+1}} [ \la \nabla f(x_{k}), y - x_{k+1} \ra + \Psi(y) - \Psi(x_{k+1})].
\ea
$$
Note that
$$
\ba{rl}
& d(y) - d(x_{k+1}) - \la \nabla d(x_k), y - x_{k+1} \ra \\
\\
\refLE{def-SConv} & d(y) - d(x_{k}) - \la \nabla d(x_k), x_{k+1} - x_{k} \ra - \half \| x_{k+1} - x_k \|^2 - \la \nabla d(x_k), y - x_{k+1} \ra\\
\\
= & r_k(y) - \half \| x_{k+1} - x_k \|^2.
\ea
$$
Thus,
$$
\ba{rl}
& r_{k+1}(y) - r_k(y) \; \leq \; {1 \over 2 L_{k+1}} [ \la \nabla f(x_{k}), y - x_{k+1} \ra + \Psi(y) - \Psi(x_{k+1})]  - \half \| x_{k+1} - x_k \|^2\\
\\
= & {1 \over 2 L_{k+1}} \Big[ \Psi(y) - \Psi(x_{k+1}) - \la \nabla f(x_{k}), x_{k+1} - x_k \ra - L_{k+1} \| x_{k+1} - x_k \|^2 \\
\\
 & + \la \nabla f(x_{k}), y - x_{k} \ra\Big]\\
\\
\leq & {1 \over 2 L_{k+1}} \Big[ \Psi(y) - \Psi(x_{k+1}) + f(x_k) - f(x_{k+1}) + \la \nabla f(x_{k}), y - x_{k} \ra\Big].
\ea
$$
Thus, we obtain inequality
$$
\ba{rcl}
{1 \over 2 L_{k+1}} \tilde f(x_{k+1}) + r_{k+1}(y) & \leq & {1 \over 2 L_{k+1}} \left( f(x_k) + \la \nabla f(x_k), y - x_k \ra + \Psi(y) + {\epsilon \over 2} \right) + r_k(y).
\ea
$$
Summing up these inequalities, we obtain
$$
\ba{rcl}
\tilde f^*_k & \leq {1 \over S_k} \sum\limits_{i=0}^k {1 \over L_{i+1}} [f(x_i) + \la \nabla f(x_i), y - x_i \ra] + \Psi(y) + {\epsilon \over 2} + {2 \over S_k} r_0(y).
\ea
$$
It remains to note that by inequality (\ref{eq-Bound2L}) we have $S_k \geq {k+1 \over \gamma_f(\epsilon)}$.
\qed

Theorem \ref{th-PGM} allows us to bound the number of steps of method (\ref{met-PGM}) in terms of the {\em complexity gauge}. Indeed, in order to get $\epsilon$-solution of our problem, we need to perform
\beq\label{eq-CGComp1}
\ba{rcl}
{4 \over \epsilon} \gamma_f(\epsilon) \beta(x_0,x^*) & = & 
{4 \beta(x_0,x^*) \over s_f^2\left({\epsilon \over 2}\right)} 
\ea
\eeq
iterations at most. As we have already discussed, this is a {\em universal} complexity bound, which works for all problem classes defined by the uniform smoothness conditions.

In the direct form (\ref{eq-UComp1}), the performance guarantee (\ref{eq-CGComp1}) looks as follows:
\beq\label{eq-CGComp2}
\ba{rcl}
2 \hat \mu \left( 2 \sqrt{D \over k} \right) & \leq & \epsilon, \quad D \Def \beta(x_0,x^*).
\ea
\eeq

Consider now the Dual Gradient Method.
\beq\label{met-DGM}
\ba{|l|}
\hline \\
\hspace{15ex} \mbox{\bf Universal Dual Gradient Method (DGM)}\\
\\
\hline \\
\mbox{{\bf Initializatin.} Choose $L_0 > 0$ and accuracy $\epsilon > 0$. Define $\phi_0(x) = \beta(x_0,x)$.}\\
\\
\mbox{{\bf For $k \geq 0$ do:}}\\
\\
\ba{rl}
1. & \mbox{Find the smallest $i_k \geq 0$ such that for the point}\\
\\
& \mbox{$x_{k,i_k} = \argmin\limits_{x \in \dom \Psi} \Big\{ \phi_k(x) + {1 \over 2^{i_k} L_k} [f(x_k) + \la \nabla f(x_k), x - x_k \ra + \Psi(x) ]\Big\}$}\\
\\
& \mbox{we have $\tilde f \left( {\cal B}_{2^{i_k}L_k}(x_{k,i_k})\right) \leq \psi^*_{2^{i_k}L_k}(x_{k,i_k}) + {\epsilon \over 2}$.}\\
\\
2. & \mbox{Set $x_{k+1} = x_{k,i_k}$, $y_k = {\cal B}_{2^{i_k}L_k}(x_{k,i_k})$, $L_{k+1} = 2^{i_k-1} L_k$, and}\\
\\
& \phi_{k+1}(x) = \phi_k(x) + {1 \over 2 L_{k+1}} [ f(x_k) + \la \nabla f(x_k), x - x_k \ra + \Psi(x)].\\
\ea\\
\\
\hline
\ea
\eeq

Denote $S_k = \sum\limits_{i=0}^k {1 \over L_{i+1}}$, $\tilde f^*_k = {1 \over S_k} \sum\limits_{i=0}^k {1 \over L_{i+1}} \tilde f(y_i)$, and $\phi^*_k = \min\limits_{y \in \dom \Psi} \phi_k(y) = \phi_k(x_k)$.
\BT\label{th-DGM}
Let $f(\cdot)$ satisfy Assumption \ref{ass-Dom} and $L_0 \leq \gamma_f(\epsilon)$. Then, for all $k \geq 0$, we have $L_k \leq \gamma_f(\epsilon)$. Moreover, for all $k\geq 0$, it holds
\beq\label{eq-RateDGM}
\ba{rcl}
\sum\limits_{i=0}^k {1 \over 2 L_{i+1}} \tilde f(y_i) & \leq & \phi^*_{k+1} + S_k \cdot {\epsilon \over 4}.
\ea
\eeq
Therefore, $\tilde f^*_k - \tilde f(x^*) \leq {\epsilon \over 2} + {2 \gamma_f(\epsilon) \over k+1} d(x_0,x^*)$.
\ET
\proof
Let us prove (\ref{eq-RateDGM}) by induction. For $k = 0$, we have
$$
\ba{rl}
& {1 \over 2L_1} \tilde f(y_0) - S_0 \cdot {\epsilon \over 4} \; = \; {1 \over 2L_1} \Big[ \tilde f(y_0) - {\epsilon \over 2} \Big] \; \leq \; {1 \over 2^{i_0}L_0} \psi^*_{2^{i_0}L_0}(y_0)\\
\\
= & {1 \over 2^{i_0}L_0} \Big[ f(x_0) + \la \nabla f(x_0), y_0 - x_0 \ra + \Psi(y_0) \Big] + \beta(x_0,y_0)\\
\\
= & \phi_1(y_0) \; = \; \min\limits_{x \in \dom \Psi} \phi_1(x) \; = \; \phi^*_1.
\ea
$$

Assume now that (\ref{eq-RateDGM}) holds for some $k \geq 0$. Note that
$$
\ba{rcl}
\phi_k(x) & \geq & \phi_k(x_k) + \beta(x_k,x) \; = \; \phi^*_k + \beta(x_k,x), \quad x \in \dom \Psi.
\ea
$$
(See, for example, Lemma 3 in \cite{UGM}.) Therefore,
$$
\ba{rl}
& \phi^*_{k+2} \; = \; \min\limits_{x \in \dom \Psi} \phi_{k+2}(x)\\
= & \min\limits_{x \in \dom \Psi} \Big\{ \phi_{k+1}(x) + {1 \over 2 L_{k+2}} [ f(x_{k+1}) + \la \nabla f(x_{k+1}), x - x_{k+1} \ra + \Psi(x) ]\Big\}\\
\\
\geq & \min\limits_{x \in \dom \Psi} \Big\{ \phi^*_{k+1} + \beta(x_{k+1},x) + {1 \over 2 L_{k+2}} [ f(x_{k+1}) + \la \nabla f(x_{k+1}), x - x_{k+1} \ra + \Psi(x) ]\Big\}\\
\\
\geq & \phi_{k+1}^* + {1 \over 2 L_{k+2}} \Big[ \tilde f(y_{k+1}) - {\epsilon \over 2} \Big] \; \refGE{eq-RateDGM} \; - S_{k+1} \cdot {\epsilon \over 4} + \sum\limits_{i=0}^{k+1} {1 \over 2 L_{i+1}} \tilde f(y_i).
\ea
$$
Thus, (\ref{eq-RateDGM}) is proved. 

Using the same arguments as in Theorem \ref{th-PGM}, we can prove that $L_k \leq \gamma_f(\epsilon)$. 
It remains to note that
$$
\ba{rcl}
\half S_k \tilde f^*_k & = & \sum\limits_{i=0}^k {1 \over 2L_{i+1}} \tilde f(y_i) \; \refLE{eq-RateDGM} \; \min\limits_{x \in \dom \Psi} \phi_{k+1}(x) + S_k \cdot {\epsilon \over 4} \\
\\
& \leq & \half S_k \tilde f(x^*) + \beta(x_0,x^*) + S_k \cdot {\epsilon \over 4}. \QR
\ea
$$

Thus, the rate of convergence of the Dual Gradient Method (\ref{met-DGM}) is the same as that of the primal one. Therefore, it has the same performance guarantees (\ref{eq-CGComp1}), (\ref{eq-CGComp2}).

\section{Fast Universal Gradient Method}\label{sc-UFast}
\SetEQ

Consider now a modified version of Universal Fast Gradient Method proposed in \cite{UGM}, where its universality was justified for functions with H\"older continuous gradients. In this section,
we are going to derive for it a direct complexity bound.
For the reasonings below, the simplifying assumption $\dom f = \E$ is still valid.

\beq\label{met-UFGM}
\ba{|l|}
\hline \\
\hspace{14ex} \mbox{\bf Universal Fast Gradient Method (UFGM)}\\
\\
\hline \\
\mbox{{\bf Initialization.} Choose $L_0 > 0$ and $\epsilon >0$.}\\
\\
\mbox{Define $\phi_0(x) = \beta(x_0,x)$, $y_0 = x_0$, and $A_0 = 0$.}\\
\\
\hline \\
\mbox{\bf For $k \geq 0$ do:}\\
\\
\mbox{{\bf 1.} Find $v_k = \argmin\limits_{x \in \dom \Psi} \phi_k(x)$.}\\
\\
\mbox{{\bf 2.} Find the smallest $i_k \geq 0$ such that coefficient $a_{k+1,i_k} > 0$, computed from}\\
\\
\mbox{equation \fbox{$a_{k+1,i_k}^2 = {1 \over 2^{i_k}L_k}(A_k + a_{k+1,i_k})$} and used in the definitions}\\
\\
A_{k+1,i_k} = A_k + a_{k+1,i_k}, \quad \tau_{k,i_k} = {a_{k+1,i_k} \over A_{k+1,i_k}}, \quad x_{k+1,i_k} = \tau_{k,i_k} v_k + (1-\tau_{k,i_k})y_k,\\
\\
\hspace{5ex} \hat x_{k+1, i_k} = \argmin\limits_{y \in \dom \Psi} \{ \beta(v_k,y) + a_{k+1,i_k} [ \la \nabla f(x_{k+1,i_k}), y + \Psi(y)] \},\\
\\
\mbox{and $ y_{k+1,i_k} = \tau_{k,i_k} \hat x_{k+1,i_k} + (1-\tau_{k,i_k}) y_k$, ensures the following relation:}\\
\\
\hspace{5ex}\ba{rl}
f(y_{k+1,i_k}) \; \leq & f(x_{k+1,i_k}) + \la \nabla f(x_{k+1,i_k}), y_{k+1,i_k} - x_{k+1,i_k} \ra\\ \\
& + 2^{i_k-1} L_k \| y_{k+1,i_k} - x_{k+1,i_k} \|^2 + {\epsilon \over 2} \tau_{k,i_k}.
\ea \\
\\
\mbox{{\bf 3.} Set $x_{k+1} = x_{k+1,i_k}$, $y_{k+1} = y_{k+1,i_k}$, $a_{k+1} = a_{k+1,i_k}$, and $\tau_{k} = \tau_{k,i_k}$.}\\
\\
\mbox{Define $A_{k+1} = A_k + a_{k+1}$, $L_{k+1} = 2^{i_k} L_k$, and}\\
\\
\hspace{5ex}\phi_{k+1}(x) = \phi_k(x) + a_{k+1} [ f(x_{k+1}) + \la \nabla f(x_{k+1}), x - x_{k+1} \ra + \Psi(x)].\\
\\
\hline
\ea
\eeq
Note that in \cite{UGM} at Step 3, we used a more aggressive rule $L_{k+1}=2^{i_k-1}L_k$.

\BT\label{th-UFGM}
Let $\dom f = \E$ and $f(\cdot)$ satisfies Assumption \ref{ass-Dom}. If $L_0 \leq 2 \gamma_f(\epsilon)$, then all iterations of method (\ref{met-UFGM}) are well defined and
\beq\label{eq-AUBound}
\ba{rcl}
L_{k+1} & \leq & 2 \gamma_f \left(\epsilon \, \tau_k \right), \quad \tau_k \; = \; {a_{k+1} \over A_{k+1}}, \quad k \geq 0.
\ea
\eeq
Moreover, for all $k \geq 0$, we have 
\beq\label{eq-RTau}
\ba{rcl}
\tau_{k+1} & \leq & \tau_k \; \leq \; {2 \over k+1},
\ea
\eeq 
and 
\beq\label{eq-UFGM}
\ba{rcl}
A_k \Big(\tilde f(y_k) - {1 \over 2} \epsilon \Big) & \leq & \phi^*_k \Def \min\limits_{x \in \dom \Psi} \phi_k(x).
\ea
\eeq
\ET
\proof
Let us prove first that the method is well defined. Denote $t_{k,i_k} = \epsilon \tau_{k,i_k}$ and $r_{k,i_k} = s_f(t_{k,i_k})$. Thus, $t_{k,i_k} = \hat \mu (r_{k,i_k})$ and in view of Lemma \ref{lm-MBig}, the termination condition of Step 2 is surely activated if 
\beq\label{eq-LKBound}
\ba{rcl}
2^{i_k} L_k & \geq & \gamma_f(t_{k,i_k}) \; \refEQ{def-GammaF} \; {t_{k,i_k} \over s_f^2\left(\half t_{k,i_k} \right)}.
\ea
\eeq
Note that for any $i \geq 0$, we have $a_{k+1,i+1} \leq a_{k+1,i}$. Therefore, $t_{k+1,i+1} \leq t_{k+1,i}$, and we conclude that
$$
\ba{rcl}
s_f\left(\half t_{k,i_k}\right) & = & s_f\left( {t_{k,i_k} \over 2 t_{k,0}} \cdot t_{k,0} \right) \; \refGE{eq-GBound} \; {t_{k,i_k} \over t_{k,0}} \cdot s_f\left( \half t_{k,0} \right).
\ea
$$
Hence, the termination condition is surely activated when
$$
\ba{rcl}
{t^2_{k,0} \over s_f^2\left( \half t_{k,0} \right)}  & \leq & 2^{i_k} L_k t_{k,i_k} \; = \; \epsilon \cdot 2^{i_k} L_k {a_{k+1,i_k} \over A_{k+1,i_k}} \; = \; {\epsilon \over a_{k+1,i_k}}.
\ea
$$
It remains to note that and $a_{k+1,i} \to 0$ as $i \to \infty$. Hence, the line-search process at Step~2 is finite.

Let us prove now the relation (\ref{eq-UFGM}). For $k=0$, it is evident since $A_0 = 0$ and $\phi_0^* = 0$. Assume that it is valid for some $k \geq 0$. Then, for any $y \in \dom \Psi$, we have
$$
\ba{rcl}
\phi_k(y) & \geq & \phi_k^* + \beta(v_{k},y) \; \refGE{eq-UFGM} \; A_k \Big(\tilde f(y_k) - \half\epsilon \Big) + \beta(v_k,y).
\ea
$$
Therefore,
$$
\ba{rcl}
\phi_{k+1}(y) + \half A_k \epsilon & \geq & \beta(v_k,y) + A_k \Big( f(x_{k+1}) + \la \nabla f(x_{k+1}), y_k - x_{k+1} \ra + \Psi(y_k) \Big)\\
\\
& & + a_{k+1} [ f(x_{k+1}) + \la \nabla f(x_{k+1}), y - x_{k+1} \ra + \Psi(y)]\\
\\
& = & \beta(v_k,y) + A_k \Big( f(x_{k+1}) + \Psi(y_k) \Big)\\
\\
& & + a_{k+1} [ f(x_{k+1}) + \la \nabla f(x_{k+1}), y - v_k \ra + \Psi(y)].
\ea
$$

In view of definition of the point $\hat x_{k+1}$, we have
$$
\ba{rcl}
\phi^*_{k+1} + \half A_k \epsilon & \geq & \beta(v_k, \hat x_{k+1}) + A_k \Big( f(x_{k+1}) + \Psi(y_k) \Big)\\
\\
& & + a_{k+1} [ f(x_{k+1}) + \la \nabla f(x_{k+1}), \hat x_{k+1} - v_k \ra + \Psi(\hat x_{k+1})]\\
\\
& \refGE{eq-BGrow} & \half \| \hat x_{k+1} - v_k \|^2 + A_{k+1} f(x_{k+1}) + A_{k+1} \Psi(y_{k+1}) \\
\\
& & + a_{k+1} \la \nabla f(x_{k+1}), \hat x_{k+1} - v_k \ra.
\ea
$$
Since $\hat x_{k+1} - v_k = {1 \over \tau_k}(y_{k+1} - x_{k+1})$, we obtain
$$
\ba{rcl}
\phi^*_{k+1} + \half A_k \epsilon & \geq & {1 \over 2 \tau_k^2} \| y_{k+1}- x_{k+1} \|^2 + 
A_{k+1} f(x_{k+1}) + A_{k+1} \Psi(y_{k+1}) \\
\\
& & + A_{k+1} \la \nabla f(x_{k+1}), y_{k+1} - x_{k+1} \ra\\
\\
& =  & A_{k+1} ( f(x_{k+1}) + \la \nabla f(x_{k+1}), y_{k+1} - x_{k+1} \ra\\
\\
& & + 2^{i_k-1}L_k \| y_{k+1} - x_{k+1} \|^2 + \Psi(y_{k+1})) \\
\\
& \geq & A_{k+1} \Big( f(y_{k+1}) - {\epsilon \over 2} \tau_k  + \Psi(y_{k+1}) \Big) \; = \; A_{k+1} \tilde f(y_{k+1}) - {\epsilon \over 2} a_{k+1}.
\ea
$$
Thus, inequality (\ref{eq-UFGM}) is proved for all $k \geq 0$.

Let us prove now the inequality (\ref{eq-AUBound}).
If $k = 0$, then $A_1 = a_1$ and $\tau_{0,i} = 1$ for all $i \geq 0$.
Hence, if $i_0 = 0$, then $L_1 = L_0$ and this inequality is valid in view of our assumption. If $i_0 \geq 1$, then in view of Lemma \ref{lm-MBig}, we have
$$
\ba{rcl}
\half L_1 & = & 2^{i_0-1} L_0 \leq \gamma_f(\epsilon).
\ea
$$
Hence, we get that (\ref{eq-AUBound}) is valid for $k = 0$.

Further, for some $A \geq 0$ and $L > 0$, consider the functions $a_L(A)$ and $\tau_L(A)$, defined by the following equations:
$$
\ba{rcl}
{a_L^2(A) \over A + a_L(A)} & = & {1 \over L}, \quad \tau_L(A) \; = \; {a_L(A) \over A + a_L(A)}.
\ea
$$
From the second equation, we get $a_{L}(A) = {\tau_L(A) A \over 1 - \tau_L(A)}$, and 
$$
\ba{rcl}
\tau^2_L(A) A^2 & = & {1 \over L} (A + a_L(A))(1-\tau_L(A))^2 \; = \; {1 \over L} A(1-\tau_L(A)).
\ea
$$
Hence, ${\tau^2_L(A) \over 1 - \tau_L(A)} = {1 \over LA}$, which implies that $\tau_L(A)$ is a decreasing function of $LA$.

Assume now that inequality (\ref{eq-AUBound}) is true for some $k \geq 0$, that is  
$$
\ba{rcl}
L_{k+1} & \refLE{eq-AUBound} & 2 \gamma_f \left(\epsilon \cdot \tau_k \right), \quad \tau_k \; = \; \tau_{L_{k+1}} (A_{k}).
\ea
$$
Let us look first at the case $i_{k+1} = 0$. Then $L_{k+2} = L_{k+1}$ and, since $A_{k+1} \geq A_k$, we have
$$
\ba{rcl}
\tau_{k+1} & = & \tau_{L_{k+2}}(A_{k+1}) \; = \;\tau_{L_{k+1}}(A_{k+1}) \; \leq \; \tau_{L_{k+1}}(A_k) \; = \; \tau_k.
\ea
$$
Therefore, in view of inequality (\ref{eq-GNon}), we have $L_{k+2} \leq 2 \gamma_f( \epsilon \cdot \tau_{k+1})$.

Consider now the case $i_{k+1} \geq 1$.
In view of Lemma \ref{lm-MBig}, we have
$$
\ba{rcl}
2^{i_{k+1}-1} L_{k+1} & \leq & \gamma_f(t_{k+1, i_{k+1}-1}) \; \leq \; \gamma_f(t_{k+1, i_{k+1}}).
\ea
$$
Hence, $L_{k+2} \leq 2 \gamma_f(\epsilon \cdot \tau_{k+1})$, and we conclude that inequality (\ref{eq-AUBound}) is valid for all $k \geq 0$.

Since $L_{k+2} \geq L_{k+1}$, as a by-product of our reasoning, we get the following relation:
$$
\ba{rcl}
\tau_{k+1} & = & \tau_{L_{k+2}}(A_{k+1}) \; \leq \; \tau_{L_{k+1}}(A_{k}) \; = \tau_k, \quad k \geq 0,
\ea
$$
which is the first inequality in (\ref{eq-RTau}).

For proving the second part of (\ref{eq-RTau}), note that
${a_{i+1} \over A_{i+1}^{1/2}} = {1 \over \sqrt{L_{i+1}}}$.
Since the univariate function $\xi(\tau) = \sqrt{\tau}$, $\tau \geq 0$, is concave,
we have
$$
\ba{rcl}
A_{i+1}^{1/2} - A_i^{1/2} & \geq & \xi'(A_{i+1})(A_{i+1} - A_i) \; = \; {A_{i+1} - A_i \over 2 A_{i+1}^{1/2}} \; = \; {1 \over 2 \sqrt{L_{i+1}}}, \quad i \geq 0.
\ea
$$
Summing up these inequalities for $k \geq 1$, we get
$$
\ba{rcl}
A_k & \geq & {1 \over 4} \Big[ \sum\limits_{i=1}^k {1 \over \sqrt{ L_i}}\Big]^2 \; \geq \; {k^2 \over 4 L_k}.
\ea
$$
This implies that $\tau_k = {1 \over \sqrt{L_{k+1} A_{k+1}}}\leq {2 \over k+1}$, which gives us the second part of (\ref{eq-RTau}).
\qed

Let us derive some consequences of Theorem \ref{th-UFGM}.
Since $\phi_{k+1}(y) \leq A_{k+1} \tilde f(y) + \beta(x_0,y)$ for all $y \in \dom \Psi$, we obtain
\beq\label{eq-DF}
\ba{rcl}
\tilde f(y_{k+1}) - \tilde f(x^*) & \refLE{eq-UFGM} & {D \over A_{k+1}} + {\epsilon \over 2}, \quad k \geq 0,
\ea
\eeq
where $D \Def \beta(x_0,x^*)$.

Now, we can get the direct complexity bounds for method (\ref{met-UFGM}).
\BT\label{th-RateUFGM}
1. If for some $k \geq 0$, we have
\beq\label{eq-AssA}
\ba{rcl}
A_{k+1} & \leq & {2 \over \epsilon}D,
\ea
\eeq
then 
\beq\label{eq-EMu}
\ba{rcl}
{2 \over \tau_k} \hat \mu \Big( 2  \tau_k^{3/2} D^{1/2}\Big) & \geq & \epsilon.
\ea
\eeq

2. Let the number of iterations $k \geq 0$ of method (\ref{met-UFGM}) satisfy inequality
\beq\label{eq-RateUFGM}
\ba{rcl}
(k+1)\,  \hat \mu \Big( 2 \left( {2 \over k+1} \right)^{3/2} D^{1/2} \Big) & \leq & \epsilon.
\ea
\eeq
Then $A_{k+1} \geq {2 \over \epsilon}D$ and $\tilde f(y_{k+1}) - \tilde f(x^*) \leq \epsilon$.
\ET
\proof
1. By the rules of Step 2 of method (\ref{met-UFGM}), we have
$$
\ba{rcl}
{1 \over A_{k+1} \tau_k^2} & = & L_{k+1} \; \refLE{eq-AUBound} \; 2 \gamma_f(\epsilon \, \tau_k) \; \refEQ{def-GammaF} \; {2 \epsilon \, \tau_k \over s_f^2\left( \half \epsilon \, \tau_k \right)}.
\ea
$$
Therefore,
${2 \tau^3_k \over s_f^2\left( \half \epsilon \, \tau_k \right)} \geq {1 \over \epsilon A_{k+1}} \; \refGE{eq-AssA} \; {1 \over 2D}$.
In other words, $2 D^{1/2} \tau_k^{3/2} \geq s_f\left( \half \epsilon \, \tau_k \right)$, which is the same as (\ref{eq-EMu}).

2. Let (\ref{eq-RateUFGM}) hold for some $k \geq 0$. Assume that the strict variant of inequality (\ref{eq-AssA}) is valid.
Note that, by inequality (\ref{eq-HMuUp}), the left-hand side of inequality (\ref{eq-EMu}) is monotonically increasing in $\tau_k$. Then, using the second part of inequality (\ref{eq-RTau}), we come to the following relation:
$$
\ba{rcl}
\epsilon & < & (k+1)\, \hat \mu \Big( 2  \left( {2 \over k+1} \right)^{3/2} D^{1/2}\Big),
\ea
$$
which contradicts to condition (\ref{eq-RateUFGM}). Hence, our assumption is wrong and $A_{k+1} \geq {2 \over \epsilon}D$. Thus, in view of (\ref{eq-DF}), point $y_{k+1}$ satisfies the desired inequality.
\qed

This theorem has several important consequences. In particular, we can estimate the total number of additional calls of oracle during the search procedure of Step 2.
Denote
$$
\ba{rcl}
N_k & = & \sum\limits_{p=0}^k i_p \; = \; \sum\limits_{p=0}^k \log_2 {L_{p+1} \over L_p} \; = \; \log_2 {L_{k+1} \over L_0}, \quad k \geq 0.
\ea
$$
Assume that $A_{k+1} \leq {2 \over \epsilon} D$. Then, in view of inequality (\ref{eq-EMu}), we have
$$
\ba{rcl}
\epsilon & \leq & {2 \over \tau_k} \hat \mu \Big( 2  \tau_k^{3/2} D^{1/2}\Big) \; \refLE{eq-HMuUp} \; 2 \tau^{1/2}_k \hat \mu( 2 D^{1/2}).
\ea
$$
Therefore, we get the following bound:
\beq\label{eq-UpN}
\ba{c}
N_k \; = \; \log_2 {L_{k+1}\over L_0} \; \refLE{eq-AUBound} \; 1+ \log_2 {\gamma_f \left(\epsilon \tau_k \right) \over L_0} \\
\\
\refLE{eq-GNon} \; 1 + \log_2 \gamma_f \left( {\epsilon^3 \over 4 \hat \mu^2 ( 2 D^{1/2})}\right) - \log_2 L_0.
\ea
\eeq

\section{Conclusion}\label{sc-Conc}
\SetEQ

Let us discuss some consequences of the main results of this paper.

{\bf 1.} Condition (\ref{eq-RateUFGM}) can be seen as a rate of convergence of method (\ref{met-UFGM}). In view of inequality (\ref{eq-HMuUp}), it admits the following sufficient condition:
\beq\label{eq-Worst}
\ba{rcl}
{2^{3/2} \over \sqrt{k+1}} \; \hat \mu(2 D^{1/2}) & \leq & \epsilon,\quad k \geq 1.
\ea
\eeq
This is a very general rate of convergence, which does not depend on the problem classes containing $f(\cdot)$. It depends only on curvature of function $f(\cdot)$ with respect to the chosen norm $\| \cdot \|$.

{\bf 2.} Let us show how the condition (\ref{eq-RateUFGM}) works for the standard problem classes. Assume that $f \in \C^{1,\nu}(\E)$ with $\nu \in [0,1]$. Then $\hat \mu(t) \refLE{eq-MUp} {L_{\nu} \over 1+\nu}t^{1 + \nu}$. Therefore, a sufficient condition for inequality~(\ref{eq-RateUFGM}) is as follows:
\beq\label{eq-SufNu}
\ba{rcl}
{L_{\nu} \over 1+\nu} \Big( 2^{5/2} D^{1/2} \Big)^{1+\nu} \left({1 \over k+1} \right)^{1 + 3 \nu \over 2}& \leq & \epsilon.
\ea
\eeq
It is known that this is the optimal rate of convergence for the corresponding problem class \cite{NY}. Our method achieves it simultaneously for all $\nu \in [0,1]$.

{\bf 3.} Comparing the bounds (\ref{eq-CGComp2}) and (\ref{eq-RateUFGM}), we can see that 
$$
\ba{rcl}
(k+1)\,  \hat \mu \Big( 2 \left( {2 \over k+1} \right)^{3/2} D^{1/2} \Big) & \refLE{eq-HMuUp}  & 2\,  \hat \mu \Big( 2 \left( {2 D\over k+1} \right)^{1/2} \Big), \quad k \geq 1.
\ea
$$
Thus, the Fast Gradient Method is usually better than the simple gradient methods, especially on good problem classes. Its superiority is maximal for the class $\C^{1,1}(\E)$, where $\hat \mu(t) \leq {L_1 \over 2} t^2$. For the class $\C^{1,0}(\E)$ with $\hat \mu(t) \leq L_0t$, they are both optimal \cite{NY}.

{\bf 4.} The advantages of our framework can be demonstrated by non-standard problem classes. For example, the complexity bounds for the problem class of Example \ref{ex-Sum}, can be readily derived by substituting the particular form of its modulus of uniform smoothness~(\ref{eq-Sum}) into the inequalities (\ref{eq-CGComp2}) and (\ref{eq-RateUFGM}).

{\bf 5.} For the future research, the most interesting direction is related to development of universal method based on the exact modulus of uniform convexity.

\subsection*{Acknowledgement}

The author is thankful to Ion Necoara and Nikita Doikov for the interesting and motivating discussions related to the topic of this paper.

We would like to thank for the support the National Research, Development and
Innovation Office (NKFIH) under grant number 2024-1.2.3-HURIZONT-2024-00030.

\subsection*{Conflict of interests}

The author declares no conflicts of interests.

\end{document}